\documentstyle[12pt]{article}
\def\Bbb R{{\rm \bf R}}
\def\proclaim#1{\vskip2mm{\bf #1}\em}
\def\endproclaim{\em \vskip2mm}
\def\tag#1{\eqno(#1)}
\def\gathered{\begin{array}{c}}
\def\endgathered{\end{array}}
\def\text{\mbox}

\begin{document}

\title {The enclosure method for inverse obstacle scattering problems
with dynamical data over a finite time interval}
\author{Masaru IKEHATA\\
Department of Mathematics,
Graduate School of Engineering\\
Gunma University, Kiryu 376-8515, JAPAN}
\date{18 September 2009}
\maketitle

\begin{abstract}
A simple method for some class of inverse obstacle scattering
problems is introduced.  The observation data are given by a wave
field measured on a known surface surrounding unknown obstacles
over a {\it finite} time interval. The wave is generated by an
initial data with compact support outside the surface. The method
yields the distance from a given point outside the surface to
obstacles and thus more than the convex hull.

\noindent
AMS: 35R30

\end{abstract}

%\tableofcontents

\section{Introduction and statements of the results}
The aim of this paper is to introduce a simple method for some
class of inverse obstacle scattering problems in which some {\it
dynamical data} over a {\it finite} time interval are used as
the observation data.

In order to explain the essence of the idea
we consider two inverse obstacle scattering problems whose
governing equations are given by the wave equations in three
dimensions.

\subsection{Sound hard obstacles}

Let $0<T<\infty$.  Let $D\subset\Bbb R^3$ be a bounded open set with smooth boundary
such that $\Bbb R^3\setminus\overline D$ is connected.
Given $f\in L^2(\Bbb R^3)$ with compact support satisfying
$\text{supp}\,f\cap\overline D=\emptyset$ let $u=u(x,t)$
satisfy the initial boundary value problem:
$$\begin{array}{c}
\displaystyle
\partial_t^2u-\triangle u=0\,\,\text{in}\,(\Bbb R^3\setminus\overline D)\times\,]0,\,T[,\\
\\
\displaystyle
\frac{\partial u}{\partial\nu}=0\,\,\text{on}\,\partial D\times\,]0,\,T[,\\
\\
\displaystyle
u(x,0)=0\,\,\text{in}\,\Bbb R^3\setminus\overline D,\\
\\
\displaystyle
\partial_tu(x,0)=f(x)\,\,\text{in}\,\Bbb R^3\setminus\overline D.
\end{array}
\tag {1.1}
$$
Here we denote the unit outward normal to $\partial D$
by the symbol $\nu$.

Let $\Omega$ be a bounded domain with smooth boundary such that $\overline D\subset\Omega$
and $\Bbb R^3\setminus\overline\Omega$ is connected.
We denote the unit outward normal to $\partial\Omega$
by $\nu$ again.  The $\partial\Omega$ is considered as the location of the receivers
of the acoustic wave produced by an emitter located at the support of $f$.
In this paper first we consider the following problem.

$\quad$

{\bf\noindent Inverse Problem I.}
Assume that $D$ is unknown.
Extract information about the location and shape of $D$
from $u$ on $\partial\Omega\times]0,\,\,T[$ for some fixed {\it known} $f$
satisfying $\text{supp}\,f\cap\overline\Omega=\emptyset$ and $T<\infty$.

$\quad$

This is a quite natural problem, however, to my best knowledge, it seems that no attempt has been done.
Clearly the main obstruction is the {\it finiteness} of $T$ and $f$ is fixed.

Note that $u$ in $(\Bbb R^3\setminus\overline\Omega)\times]0,\,\,T[$ can be computed from
$u$ on $\partial\Omega\times\,]0,\,T[$ by
the formula
$$\displaystyle
u=z\,\,\text{in}\,(\Bbb R^3\setminus\overline\Omega)\times\,]0,\,T[
\tag {1.2}
$$
where $z$ solves the initial boundary value problem in $\Bbb R^3\setminus\overline\Omega$:
$$\begin{array}{c}
\displaystyle
\partial_t^2z-\triangle z=0\,\,\text{in}\,(\Bbb R^3\setminus\overline\Omega)\times\,]0,\,T[,\\
\\
\displaystyle
z=u\,\,\text{on}\,\partial\Omega\times\,]0,\,T[,\\
\\
\displaystyle
z(x,0)=0\,\,\text{in}\,\Bbb R^3\setminus\overline\Omega,\\
\\
\displaystyle
\partial_tz(x,0)=f(x)\,\,\text{in}\,\Bbb R^3\setminus\overline\Omega.
\end{array}
\tag {1.3}
$$

Thus the problem can be reformulated as

$\quad$

{\bf\noindent Inverse Problem I'.}
Extract information about the location and shape of $D$ from
$u$ in $(\Bbb R^3\setminus\overline\Omega)\times]0,\,T[$ for some known $f$ satisfying $\text{supp}\,f\cap\overline\Omega=\emptyset$
and $T<\infty$.

$\quad$

Now we state the result.
Let $B$ be an open ball with $\overline B\cap\overline\Omega=\emptyset$.
Choose the initial data $f\in L^2(\Bbb R^3)$ in such a way that:

$\quad$

\noindent
(I1)  $f(x)=0$ a.e. $x\in\Bbb R^3\setminus B$,

$\quad$

\noindent
(I2)  there exists a positive constant $C$ such that $f(x)\ge C$ a.e. $x\in B$
or $-f(x)\ge C$ a.e. $x\in B$.

$\quad$

Let $\tau>0$ and $v\in H^1(\Bbb R^3)$ be the weak solution of
$$\displaystyle
(\triangle-\tau^2)v+f(x)=0\,\,\text{in}\,\Bbb R^3.
\tag {1.4}
$$
This means that $v$ satisfies
$$
\displaystyle
-\int_{\Bbb R^3}\nabla v\cdot\nabla\varphi dx-\tau^2\int_{\Bbb R^3}v\varphi dx
=-\int_{\Bbb R^3}f\varphi dx,\,\,\forall\varphi\in H^1(\Bbb R^3).
\tag {1.5}
$$
The $v$ is unique and is given by the explicit form
$$\displaystyle
v(x;\tau)=\frac{1}{4\pi}\int_B
\frac{e^{-\tau \vert x-y\vert}}{\vert x-y\vert}f(y)dy,\,\,x\in\Bbb R^3.
\tag {1.6}
$$

In this paper given two sets $A$ and $B$ we denote by $\text{dist}\,(A,B)$
the distance between $A$ and $B$:
$$\displaystyle
\text{dist}\,(A,B)=\inf\{\vert x-y\vert\,\vert\,x\in A, y\in B\}.
$$
If $A$ or $B$ consists of a single point, say $B=\{p\}$, we write
$\text{dist}\,(A,B)=d_A(p)$.

Set
$$\displaystyle
w(x;\tau)=\int_0^T e^{-\tau t}u(x,t)dt,\,\,x\in\Bbb R^3\setminus\overline\Omega,\,\,\tau>0.
\tag {1.7}
$$

Our result is the following extraction formula from
$w$ and $\partial w/\partial\nu$ on $\partial\Omega\times\,]0\,\,T[$ which
can be computed from the data $u$ in $(\Bbb R^3\setminus\overline\Omega)\times\,]0,\,T[$.

\proclaim{\noindent Theorem 1.1.}
If the observation time $T$ satisfies
$$
T>2\text{dist}\,(D,B)-\text{dist}\,(\Omega,B),
\tag {1.8}
$$
then there exists a $\tau_0>0$ such that, for all $\tau\ge\tau_0$
$$\displaystyle
\int_{\partial\Omega}\left(\frac{\partial v}{\partial\nu}w-\frac{\partial w}{\partial\nu}v\right)dS>0
$$
and the formula
$$\displaystyle
\lim_{\tau\longrightarrow\infty}
\frac{1}{2\tau}\log
\int_{\partial\Omega}\left(\frac{\partial v}{\partial\nu}w-\frac{\partial w}{\partial\nu}v\right)dS
=-\text{dist}\,(D,B),
\tag {1.9}
$$
is valid.
\endproclaim

\noindent
Since $\text{dist}\,(D,B)+\sqrt{\vert\partial B\vert/4\pi}$ coincides with the distance from the center of $B$ to $D$,
(1.9) yields the information about $d_D(p)$ for a given point $p$ in $\Bbb R^3\setminus\overline\Omega$.
Therefore one can extract more than the convex hull of $D$.
Note that we do not assume the special form of $f$ except for the conditions (I1) and (I2).

The restriction (1.8) on the observation time $T$ is reasonable.
Define the quantity
$$\displaystyle
l(\partial B,\partial D,\partial\Omega)
=\inf\,\{\vert x-y\vert+\vert y-z\vert\,\vert\,x\in\partial B\,, y\in\partial D,\, z\in\partial\Omega\}.
$$
This is the minimum length of the broken paths that start at $x\in\partial B$ and reflect at $y\in\partial D$ and return
to $z\in\partial\Omega$.  We have

\proclaim{\noindent Proposition 1.1.}
$$\displaystyle
2\text{dist}\,(D,B)-\text{dist}\,(\Omega,B)\ge l(\partial B,\partial D,\partial\Omega).
$$
\endproclaim
{\it\noindent Proof.}
One can find $x_0\in\partial B$ and $y_0\in\partial D$ such that $\vert x_0-y_0\vert=\text{dist}\,(D,B)$.
Let $l(x_0,y_0)=\{tx_0+(1-t)y_0\,\vert\,0<t<1\}$.  We see that $l(x_0,y_0)\cap\partial\Omega\not=\emptyset$.
Let $z_0\in l(x_0,y_0)\cap\partial\Omega$.  We have $\vert x_0-z_0\vert\ge\text{dist}\,(\Omega,B)$.
Thus $2\text{dist}\,(D,B)-\text{dist}\,(\Omega,B)\ge 2\vert x_0-y_0\vert-\vert x_0-z_0\vert$.  Since
$\vert x_0-z_0\vert=\vert x_0-y_0\vert-\vert y_0-z_0\vert$ we have
$2\text{dist}\,(D,B)-\text{dist}\,(\Omega,B)\ge\vert x_0-y_0\vert+\vert y_0-z_0\vert$.

\noindent
$\Box$

Therefore (1.8) ensures that $T>l(\partial B,\partial D,\partial\Omega)$.
This means that $T$ is greater than the {\it first arrival time}
of a signal with the unit propagation speed that starts at a point on $\partial B$ at $t=0$,
reflects at a point on $\partial D$ and goes to a point on $\partial\Omega$.
However, curiously enough in the proof of Theorem 1.1 we never make use of the {\it finite propagation property}
of the signal governed by the wave equation.

The procedure of extracting information about the location of $D$ is extremely simple and summarized as follows.

\noindent
(i)  Give an open ball $B$ with $\overline B\cap\overline\Omega=\emptyset$.
Using the initial data $f$ satisfying (I1) and (I2) generate the wave field $u$.

\noindent
(ii) Choose a large $T$, say, such that $T>2\sup\,\{\vert y-x\vert\,\vert\,y\in\Omega,x\in B\}
-\text{dist}\,(\Omega,B)$
and measure $u$ on $\partial\Omega$ over the time interval $]0,\,T[$.

\noindent
(iii) Compute the values of $z$ in a neighbourhood of $\partial\Omega$ relative to
$\Bbb R^3\setminus\overline\Omega$ over $]0,\,T[$ by solving (1.3).

\noindent
(iv) Choose a large $\tau$ and compute $w$ and $\partial w/\partial\nu$ on $\partial\Omega$ 
over the time interval $]0,\,T[$ via (1.2) and (1.7).

\noindent
(v)  Compute $v$ and $\partial v/\partial\nu$ on $\partial\Omega$ via (1.6).

\noindent
(vi)  Compute the quantity
$$\displaystyle
\frac{1}{2\tau}
\log
\int_{\partial\Omega}
\left(\frac{\partial v}{\partial\nu}w-\frac{\partial w}{\partial\nu}v\right)dS
$$
as an approximation of $-\text{dist}\,(D,B)$.

One choice of $f$ gives one information about $D$ by the procedure
(i) to (vi).
This means that we don't need to use {\it many} $f$s
to get $d_D(p)$ for a {\it single} $p$.
This is the decisive character of our procedure.

\subsection{Penetrable obstacles}

Given $f\in L^2(\Bbb R^3)$ with compact support
let $u=u(x,t)$
satisfy the initial value problem:
$$\begin{array}{c}
\displaystyle
\partial_t^2u-\nabla\cdot\gamma\nabla u=0\,\,\text{in}\,\Bbb R^3\times\,]0,\,T[,\\
\\
\displaystyle
u(x,0)=0\,\,\text{in}\,\Bbb R^3,\\
\\
\displaystyle
\partial_tu(x,0)=f(x)\,\,\text{in}\,\Bbb R^3,
\end{array}
\tag {1.10}
$$
where $\gamma=\gamma(x)=(\gamma_{ij}(x))$ satisfies

\noindent
$\bullet$  for each $i,j=1,2,3$ $\gamma_{ij}(x)=\gamma_{ji}(x)\in L^{\infty}(\Bbb R^3)$;

\noindent
$\bullet$  there exists a positive constant $C$ such that $\gamma(x)\xi\cdot\xi\ge C\vert\xi\vert^2$ for all $\xi\in\Bbb R^3$
and a. e. $x\in\Bbb R^3$.

This subsection is concerned with the extraction of information
about {\it discontinuity} of $\gamma$ from $u$ on
$\partial\Omega\times ]0,\,T[$ for some $f$ for a fixed
$T<\infty$.  However, we do not consider the completely general
case. Instead we assume:

\noindent
$\bullet$  there exists a bounded open set $D$ with a smooth boundary
such that $\gamma(x)$ a.e. $x\in\Bbb R^3\setminus D$ coincides with the $3\times 3$ identity matrix $I_3$.

Write $h(x)=\gamma(x)-I_3$ a.e. $x\in D$.
Our second inverse problem is the following.

$\quad$

{\bf\noindent Inverse Problem II.}  Assume that both $D$ and $h$ are {\it unknown}
and that one of the following two conditions is satisfied:

\noindent
(A1)  there exists a positive constant $C$ such that $-h(x)\xi\cdot\xi\ge\vert\xi\vert^2$
for all $\xi\in\Bbb R^3$ and a.e. $x\in D$;

\noindent
(A2)  there exists a positive constant $C$ such that $h(x)\xi\cdot\xi\ge\vert\xi\vert^2$
for all $\xi\in\Bbb R^3$ and a.e. $x\in D$.

\noindent
Extract information about the location and shape of $D$
from $u$ on $\partial\Omega\times]0,\,\,T[$ for some fixed {\it known} $f$
satisfying $\text{supp}\,f\cap\overline\Omega=\emptyset$ and $T<\infty$.

$\quad$

Note that $u$ in $(\Bbb R^3\setminus\overline\Omega)\times]0,\,\,T[$ can be computed from
$u$ on $\partial\Omega\times]0,\,T[$ by the exactly same formula as (1.2)
and thus the problem can be reformulated again as

$\quad$

{\bf\noindent Inverse Problem II'.}
Extract information about the location and shape of $D$ from
$u$ in $(\Bbb R^3\setminus\overline\Omega)\times\,]0,\,T[$ for some known $f$ satisfying
$\text{supp}\,f\cap\overline\Omega=\emptyset$ and $T<\infty$.

$\quad$

Now we state our second result.

\proclaim{\noindent Theorem 1.2.}
Assume that $\gamma$ satisfies (A1) or (A2).
Let $f$ satisfy (I1) and (I2) and $v$ be the weak solution of (1.4).
Let $T$ satisfies (1.8) and $w$ be given by (1.7) with solution of (1.10).
If (A1) is satisfied, then 
there exists a $\tau_0>0$ such that, for all $\tau\ge\tau_0$
$$\displaystyle
\int_{\partial\Omega}\left(\frac{\partial v}{\partial\nu}w-\frac{\partial w}{\partial\nu}v\right)dS>0
$$
and the formula
$$\displaystyle
\lim_{\tau\longrightarrow\infty}
\frac{1}{2\tau}\log
\int_{\partial\Omega}\left(\frac{\partial v}{\partial\nu}w-\frac{\partial w}{\partial\nu}v\right)dS
=-\text{dist}\,(D,B),
$$
is valid; if (A2) is satisfied, then there exists a $\tau_0>0$ such that, for all $\tau\ge\tau_0$
$$\displaystyle
\int_{\partial\Omega}\left(\frac{\partial v}{\partial\nu}w-\frac{\partial w}{\partial\nu}v\right)dS<0
$$
and the formula
$$\displaystyle
\lim_{\tau\longrightarrow\infty}
\frac{1}{2\tau}\log
\left(-\int_{\partial\Omega}\left(\frac{\partial v}{\partial\nu}w-\frac{\partial w}{\partial\nu}v\right)dS\right)
=-\text{dist}\,(D,B),
$$
is valid.
\endproclaim

Isakov \cite{ISA}
considered an inverse problem for the equation
$\partial_t^2u-\nabla\cdot\gamma\nabla u=0$ in
$\Omega\times\,]-\infty,\,T[$ with the zero initial data $u=0$
when $t<0$ and the lateral Neumann data $\partial u/\partial\nu=h$
on $\partial\Omega\times\,]-\infty,\,T[$, where $\Omega$ is the
half-space $x_3<0$; $\gamma$ takes the value $1$ for
$x\in\Omega\setminus D$ and a positive constant $k$ for $x\in D$.
The $D$ is given by a Lipschitz continuous function $d$ on $\Bbb
R^2$ with $d<0$ as $D=\{x\vert\,x_3<d(x_1,d_2)\}$.

His problem is to recover $\partial D$ by measuring $u$ on
$\Gamma\times]0,\,T[$ for an arbitrary fixed nonempty open set
$\Gamma\subset\partial\Omega$ and known $h$.  Choosing $h$ as the
lateral Neumann data of a special solution of the wave equation in
the half-space, he showed that if $k<1$, then the data $u$ on
$\Gamma\times\,]0,\,T[$ uniquely determines the part of $D$ within
the set of all points $x=(x_1,x_2,x_3)$ with $x_3>-T/2$ and
$(x_1,x_2)\in\Gamma$.
The condition $k<1$ corresponds to (A1).
It is an open problem whether or not the same conclusion holds in
the case $k>1$ which corresponds to (A2). See also \cite{ISA2}
for this point.

Rakesh \cite{R} considered an inverse problem for the equation
$\gamma\partial_t^2u-\nabla\cdot(\gamma\nabla u)=0$ in $\Bbb R^3\times]-\infty, T[$
with $u(x,t)=(t-x_3)^2_{+}$ for $t<<0$.  Here $(s)^2_{+}=s^2$ if $s>0$ and $(s)^2_{+}=0$ if $s\le 0$.
The $\gamma$ takes $1$ outside a bounded domain $D$ with smooth boundary and a positive constant $k(\not=1)$ on $D$.
Thus the governing equation has a same constant speed inside and outside $D$.
The data in his problem is the values of $u$ on $\partial\Omega\times]-\infty, T[$, where
$\Omega$ is a bounded open set of $\Bbb R^3$ with smooth boundary and satisfies $\overline D\subset\Omega$.
He showed that; if $T>6\,\text{diam}\,(\Omega)+\inf_{x\in D}x_3$ and $D$ is {\it strictly convex},
then the data $\partial u/\partial\nu$ on $\partial\Omega\times\,]-\infty, T[$ uniquely determine $D$ itself.

Isakov employs a contradiction argument and his method starts with the
uniqueness of the continuation of the solution of the wave
equation and derives an orthogonality relation that was deduced by
denying the conclusion.

Rakesh's argument is also a contradiction argument and makes use of the uniqueness of the continuation of the solution
of the wave equation.  However, the main point is an analysis of the wave front set of $u$.
See also \cite{R2} for other results.

Unlike them we do not make use of the continuation of a wave field
nor propagation of singularities argument. The method can be
considered as an application of the {\it enclosure method} which
was originally introduced for elliptic equations in \cite{IE,I1}.
Recently in \cite{I4} the author found its application to inverse
initial boundary value problems in one-space dimensional case for
the heat and wave equations. In \cite{IK2} we extended this method
to the heat equation in two and three-space dimensional cases.
Therein the initial data is zero and a {\it special heat flux}
depending on a large parameter is used.

\subsection{Further remarks and construction of the paper}

Finally we comment on some results in the context of the
Lax-Phillips scattering theory.  Lax-Phillips in \cite{LP}
established a relation between the support function of an obstacle
and the right end point of the {\it support} of the {\it scattering kernel} which is the observation data
in their theory.
Since
the support function gives the signed distance from the origin of
coordinates to the support plane of the obstacle, the result means
that one can get an estimation of the {\it convex hull} from the data.
Note that the scattering kernel is written
by using the scattered wave over the {\it infinite} time interval
that is produced by a singular plane wave at $t<<0$ far a way from
the obstacle, and thus the data is completely different from ours.

Majda\cite{M} considered the singularity of the scattering kernel
and clarified a relation between the support function and the right end point of the {\it singular support}
of the {\it back scattering kernel}.
In \cite{MT} a similar result for an obstacle with
a finite refractive index is given. The governing equation has the form
$\alpha(x)\partial_t^2u-\triangle u=0$ and $\alpha$ has a
discontinuity across the boundary of the obstacle and takes $1$
outside the obstacle.
For other results including the Maxwell equations, hyperbolic
systems, etc. we refer the reader to \cite{M, MT} and references therein.

A brief outline of this paper is as follows.
Theorems 1.1 and 1.2 are proved in Subsections 2.2 and 3.2, respectively.
The key point of the proofs is to derive a lower estimate of the integral
$$\displaystyle
\int_{\partial\Omega}\left(\frac{\partial v}{\partial\nu}w-\frac{\partial w}{\partial\nu}v\right)dS,
\tag {1.11}
$$
where $v$ is the weak solution of (1.4).
To establish the estimate we require some integral identities; these identities are found
in Subsections 2.1 and 3.1.
Using the identities, we show that, if $T$ satisfies (1.8), then the dominant part in the lower
estimate of (1.11) in Theorem 1.1 is essentially given by the integral of the square of $v$ over $D$.
We show that this last integral is comparable
with $e^{-2\tau\text{dist}\,(D,B)}$ ignoring a multiplication of a
power of $\tau$.  This is stated in Subsection 2.2 and proved in Subsection 4.1.
Note that in the proof of Theorem 1.2
instead of $v$ the integral of $\vert\nabla v\vert^2$ over $D$
plays the same role and the corresponding estimate is stated in Subsection 3.2 and proved in Subsection 4.2.
In the final section we give a conclusion of this paper and comments on further problems.

\section{The enclosure method for sound hard obstacles}

First we specify what we mean by the solution of (1.1).
We follow the notion of the weak solution described on pp. 552-566 in \cite{DL}
and use the notation therein.

By Theorem 1 on p.558 in \cite{DL}, given $u^0\in H^1(\Bbb R^3\setminus\overline D)$
and $u^1\in L^2(\Bbb R^3\setminus\overline D)$
we know that there exists a unique $u$ satisfying
$$\displaystyle
u\in L^2(0,\,T;H^1(\Bbb R^3\setminus\overline D)),\,
u'\in L^2(0,\,T;H^1(\Bbb R^3\setminus\overline D)),\,
u''\in L^2(0,\,T;(H^1(\Bbb R^3\setminus\overline D))')
$$
such that, for all $\phi\in H^1(\Bbb R^3\setminus\overline D)$
$$\displaystyle
<u''(t),\phi>
+\int_{\Bbb R^3\setminus\overline D}\nabla u(x,t)\cdot\nabla\phi(x)dx=0\,\,\text{a.e.}\,t\in\,]0,\,T[
$$
and $\displaystyle u(x,0)=u^0,\,\,u'(x,0)=u^1$.  In this section we say that this $u$
for $u^0=0$ and $u^1=f$ is the solution of (1.1).

\subsection{A basic identity}

Let $u$ be the solution of (1.1).
Define
$$\displaystyle
w(x;\tau)=\int_0^T e^{-\tau t}u(x,t)dt,\,\,\,x\in\Bbb R^3\setminus\overline D.
$$
This $w$ belongs to $H^1(\Bbb R^3\setminus\overline D)$.

Using integration by parts formula (Proposition 2 on p.558 in [DL]), we see that,
for all $\phi\in H^1(\Bbb R^3\setminus\overline D)$ the $w$ satisfies the equation
$$\begin{array}{c}
\displaystyle
\int_{\Bbb R^3\setminus\overline D}\nabla w\cdot\nabla\phi dx
+\int_{\Bbb R^3\setminus\overline D}(\tau^2w-f)\phi dx\\
\\
\displaystyle
=-e^{-\tau T}
\int_{\Bbb R^3\setminus\overline D}(u'(x,T)+\tau u(x,T))\phi dx.
\end{array}
\tag {2.1}
$$
This means that, in a weak sense $w$ satisfies
$$\begin{array}{c}
\displaystyle
(\triangle-\tau^2)w+f(x)=e^{-\tau T}
(u'(x,T)+\tau u(x,T))\,\,\text{in}\,\Bbb R^3\setminus\overline D,\\
\\
\displaystyle
\frac{\partial w}{\partial\nu}=0\,\,\text{on}\,\partial D.
\end{array}
$$

From (2.1) for $v\in C^{\infty}_0(\Bbb R^3\setminus\overline D)$, we have:$(\triangle-\tau^2)w+f(x)=e^{-\tau T}
(u'(x,T)+\tau u(x,T))$ in $\Bbb R^3\setminus\overline D$ in the sense of distribution
and hence $\triangle w\in L^2(\Bbb R^3\setminus\overline D)$.
This yields that $w\in H^2_{\text{loc}}(\Bbb R^3\setminus\overline D)$ and $(\triangle-\tau^2)w+f(x)=e^{-\tau T}
(u'(x,T)+\tau u(x,T))$ a.e. $x\in\Bbb R^3\setminus\overline D$.

Now define $\partial w/\partial\nu\vert_{\partial\Omega}$ as $\nabla w\vert_{\partial\Omega}\cdot\nu\in H^{1/2}(\partial\Omega)$,
where $\nabla w\vert_{\partial\Omega}$
is the trace of $\nabla w$ onto $\partial\Omega$.

Let $v\in H^1(\Bbb R^3)$ be the weak solution of (1.4).
For this $v$ by the same reason as above
we have $v\in H^2_{\text{loc}}(\Bbb R^3)$
and thus $\partial v/\partial\nu\vert_{\partial\Omega}\equiv
\nabla v\vert_{\partial\Omega}\cdot\nu\in H^{1/2}(\partial\Omega)$.

In this subsection we derive the following identity.
\proclaim{\noindent Proposition 2.1.}
Let $v\in H^1(\Bbb R^3)$ be the weak solution of (1.4).
It holds that
$$\begin{array}{c}
\displaystyle
\int_{\partial\Omega}\left(\frac{\partial v}{\partial\nu}w-\frac{\partial w}{\partial\nu}v\right)dS\\
\\
\displaystyle
=\int_D\vert\nabla v\vert^2dx+\tau^2\int_D\vert v\vert^2dx
+\int_{\Bbb R^3\setminus\overline D}\vert\nabla(w-v)\vert^2dx
+\tau^2\int_{\Bbb R^3\setminus\overline D}\vert w-v\vert^2dx
\\
\\
\displaystyle
+e^{-\tau T}\int_{\Bbb R^3\setminus\overline D}(w-v)(u'(x,T)+\tau u(x,T))dx
-e^{-\tau T}\int_{\Omega\setminus\overline D}
(u'(x,T)+\tau u(x,T))vdx\\
\\
\displaystyle
-\int_{\Omega\setminus\overline D}f(w-v)dx-\int_Dfvdx.
\end{array}
\tag {2.2}
$$

\endproclaim

{\it\noindent Proof.}
First we prove that
$$\begin{array}{c}
\displaystyle
\int_{\partial\Omega}\left(\frac{\partial v}{\partial\nu}w-\frac{\partial w}{\partial\nu}v\right)dS
\\
\\
\displaystyle
=\int_{\partial D}\frac{\partial v}{\partial\nu}w-\int_{\Omega\setminus\overline D}f(w-v)dx
-e^{-\tau T}\int_{\Omega\setminus\overline D}(u'(x,T)+\tau u(x,T))v(x)dx.
\end{array}
\tag {2.3}
$$
Let $\varphi\in H^1(\Omega\setminus\overline D)$ satisfy $\varphi=0$ on $\partial\Omega$ in the sense of the trace.
Since the zero extension of this $\varphi$ belongs to $H^1(\Bbb R^3\setminus\overline D)$, it follows from (2.1)
that
$$\displaystyle
\int_{\Omega\setminus\overline D}\nabla w\cdot\nabla\varphi
+\int_{\Omega\setminus\overline D}\{\tau^2 w-f+e^{-\tau T}(u'(x,T)+\tau u(x,T))\}\varphi dx=0.
\tag {2.4}
$$

Choose $\chi\in C_0^{\infty}(\Bbb R^3)$ such that $\chi(x)\equiv 1$ in a neighbourhood of $\partial\Omega$
and $\chi(x)\equiv 0$ in a neighbourhood of $\overline D$.
Since $(1-\chi)v\vert_{\Omega\setminus\overline D}$ vanishes in a neighbourhood
of $\partial\Omega$, it follows that (2.4) is valid for $\varphi=(1-\chi)v\vert_{\Omega\setminus\overline D}$.

On the other hand, since $\chi v$ vanishes in a neighbourhood of $\overline D$ and
$w\in H^2_{\text{loc}}(\Bbb R^3\setminus\overline D)$, integration by parts yields
$$\begin{array}{c}
\displaystyle
\int_{\Omega\setminus\overline D}\nabla w\cdot\nabla(\chi v)dx\\
\\
\displaystyle
=\int_{\partial\Omega}\frac{\partial w}{\partial\nu}\chi vdS
-\int_{\Omega\setminus\overline D}(\triangle w)\chi vdx\\
\\
\displaystyle
=\int_{\partial\Omega}\frac{\partial w}{\partial\nu} vdS
-\int_{\Omega\setminus\overline D}\{\tau^2 w-f+e^{-\tau T}(u'(x,T)+\tau u(x,T))\}\chi vdx.
\end{array}
$$
From this and (2.4) for $\varphi=(1-\chi)v\vert_{\Omega\setminus\overline D}$
we obtain
$$\displaystyle
\int_{\partial\Omega}\frac{\partial w}{\partial\nu}vdS
=
\int_{\Omega\setminus\overline D}\nabla w\cdot\nabla vdx
+\int_{\Omega\setminus\overline D}\{\tau^2 w-f+e^{-\tau T}(u'(x,T)+\tau u(x,T))\}v dx.
\tag {2.5}
$$

Choose $\tilde{w}\in H^1(\Omega)$ such that $\tilde{w}=w$ in $\Omega\setminus\overline D$.
Note that $\tilde{w}=w$ on $\partial\Omega$ and $\partial D$ in the sense of the trace.
Since $v\in H^2(\Omega)$, we have
$$\begin{array}{c}
\displaystyle
\int_{\partial\Omega}\frac{\partial v}{\partial\nu} wdS
=\int_{\partial\Omega}\frac{\partial v}{\partial\nu} \tilde{w}dS\\
\\
\displaystyle
=\int_{\Omega}\triangle v \tilde{w}dx+\int_{\Omega}\nabla v\cdot\nabla\tilde{w}dx\\
\\
\displaystyle
=\int_{\Omega\setminus\overline D}(\tau^2 v-f)wdx+\int_{\Omega\setminus\overline D}\nabla v\cdot\nabla wdx\\
\\
\displaystyle
+\int_{D}(\tau^2 v-f)\tilde{w}dx+\int_{D}\nabla v\cdot\nabla\tilde{w}dx.
\end{array}
$$
On the other hand we have
$$\begin{array}{c}
\displaystyle
\int_D\nabla v\cdot\nabla\tilde{w}dx
=\int_{\partial D}
\frac{\partial v}{\partial\nu}\tilde{w}-\int_D(\triangle v)\tilde{w}dx
\\
\\
\displaystyle
=\int_{\partial D}
\frac{\partial v}{\partial\nu}wdS-\int_D(\tau^2 v-f)\tilde{w}dx,
\end{array}
$$
that is
$$\displaystyle
\int_{\partial D}\frac{\partial v}{\partial\nu}wdS
=\int_D(\tau^2 v-f)\tilde{w}dx
+\int_D\nabla v\cdot\nabla\tilde{w}dx.
$$
Therefore we obtain
$$\displaystyle
\int_{\partial\Omega}\frac{\partial v}{\partial\nu} wdS
=\int_{\Omega\setminus\overline D}(\tau^2 v-f)wdx+\int_{\Omega\setminus\overline D}\nabla v\cdot\nabla wdx
+\int_{\partial D}\frac{\partial v}{\partial\nu}wdS.
\tag {2.6}
$$
A combination of (2.5) and (2.6) gives (2.3).

Write
$$\displaystyle
\int_{\partial D}\frac{\partial v}{\partial\nu}wdS
=\int_{\partial D}\frac{\partial v}{\partial\nu}(w-v)dS+\int_{\partial D}\frac{\partial v}{\partial\nu}vdS.
\tag {2.7}
$$
It follows from (1.5) and the trace theorem that, for all $\phi\in H^1(\Bbb R^3\setminus\overline D)$
$$\displaystyle
\int_{\partial D}\frac{\partial v}{\partial\nu}\phi dS
+\int_{\Bbb R^3\setminus\overline D}\nabla v\cdot\nabla\phi dx
+\tau^2\int_{\Bbb R^3\setminus\overline D}v\phi dx=\int_{\Bbb R^3\setminus\overline D}f\phi dx.
$$
Combining this with (2.1), we obtain
$$\begin{array}{c}
\displaystyle
\int_{\partial D}\frac{\partial v}{\partial\nu}\phi dS
-\int_{\Bbb R^3\setminus\overline D}\nabla(w-v)\cdot\nabla\phi dx
-\tau^2\int_{\Bbb R^3\setminus\overline D}(w-v)\phi dx\\
\\
\displaystyle
=e^{-\tau T}\int_{\Bbb R^3\setminus\overline D}(u'(x,T)+\tau u(x,T))\phi dx.
\end{array}
\tag {2.8}
$$
This means that $w-v$ satisfies, in a weak sense
$$\begin{array}{c}
\displaystyle
(\triangle-\tau^2)(w-v)=e^{-\tau T}(u'(x,T)+\tau u(x,T))\,\,\text{in}\,\Bbb R^3\setminus\overline D,\\
\\
\displaystyle
\frac{\partial}{\partial\nu}(w-v)=-\frac{\partial v}{\partial\nu}\,\,\text{on}\,\partial D.
\end{array}
$$
Substituting $w-v$ for $\phi$ in (2.8), we obtain
$$\begin{array}{c}
\displaystyle
\int_{\partial D}\frac{\partial v}{\partial\nu}(w-v)dS
=\int_{\Bbb R^3\setminus\overline D}\vert\nabla(w-v)\vert^2 dx
+\tau^2\int_{\Bbb R^3\setminus\overline D}\vert w-v\vert^2 dx\\
\\
\displaystyle
+e^{-\tau T}\int_{\Bbb R^3\setminus\overline D}(u'(x,T)+\tau u(x,T))(w-v).
\end{array}
\tag {2.9}
$$
Now from this together with (2.3), (2.7) and the identity
$$\begin{array}{c}
\displaystyle
\int_{\partial D}\frac{\partial v}{\partial\nu}vdS
=\tau^2\int_D\vert v\vert^2dx+\int_D\vert\nabla v\vert^2dx-\int_Dfvdx,
\end{array}
\tag {2.10}
$$
we obtain (2.2).

\noindent
$\Box$

In particular, choose $f$ in such a way that $\text{supp}\,f\cap\overline\Omega=\emptyset$.
Then (2.2) becomes
$$\begin{array}{c}
\displaystyle
\int_{\partial\Omega}\left(\frac{\partial v}{\partial\nu}w-\frac{\partial w}{\partial\nu}v\right)dS\\
\\
\displaystyle
=\int_D\vert\nabla v\vert^2dx+\tau^2\int_D\vert v\vert^2dx
+\int_{\Bbb R^3\setminus\overline D}\vert\nabla(w-v)\vert^2dx
+\tau^2\int_{\Bbb R^3\setminus\overline D}\vert w-v\vert^2dx
\\
\\
\displaystyle
+e^{-\tau T}\int_{\Bbb R^3\setminus\overline D}(w-v)(u'(x,T)+\tau u(x,T))dx
-e^{-\tau T}\int_{\Omega\setminus\overline D}
(u'(x,T)+\tau u(x,T))vdx.
\end{array}
\tag {2.11}
$$
This is the basic identity for the sound-hard obstacles.

\subsection{Proof of Theorem 1.1.}
Using the identity
$$\begin{array}{c}
\displaystyle
\tau^2\vert w-v\vert^2+e^{-\tau t}(w-v)(u'(x,T)+\tau u(x,T))\\
\\
\displaystyle
=\left\vert\tau(w-v)+\frac{e^{-\tau T}}{2\tau}(u'(x,T)+\tau u(x,T))\right\vert^2
-\frac{e^{-2\tau T}}{4\tau^2}\vert u'(x,T)+\tau u(x,T)\vert^2,
\end{array}
$$
we have from (2.11)
$$\begin{array}{c}
\displaystyle
\int_{\partial\Omega}\left(\frac{\partial v}{\partial\nu}w-\frac{\partial w}{\partial\nu}v\right)dS
\ge
\tau^2\int_D\vert v\vert^2dx\\
\\
\displaystyle
-\frac{e^{-2\tau T}}{4\tau^2}\int_{\Bbb R^3\setminus\overline D}\vert u'(x,T)+\tau u(x,T)\vert^2dx
-e^{-\tau T}\int_{\Omega\setminus\overline D}
(u'(x,T)+\tau u(x,T))vdx.
\end{array}
\tag {2.12}
$$

From (1.6) we have,
for all $x\in\Bbb R^3\setminus\overline B$
$$
\displaystyle
\vert v(x)\vert\le\frac{e^{-\tau d_B(x)}}{4\pi d_B(x)}\Vert f\Vert_{L^1(B)},\,\,
\vert\nabla v(x)\vert\le\frac{e^{-\tau d_B(x)}}{4\pi}\left(\tau+\frac{1}{d_B(x)^2}\right)\Vert f\Vert_{L^1(B)}.
\tag {2.13}
$$
This gives $\Vert v\Vert_{L^2(\Omega\setminus\overline D)}
=O(e^{-\tau\text{dist}\,(\Omega,B)})$ and note that
$$\displaystyle
\int_{\Bbb R^3\setminus\overline D}\vert u'(x,T)+\tau u(x,T)\vert^2dx
=O(\tau^2).
$$
Using these estimates, we obtain
$$\begin{array}{c}
\displaystyle
\int_{\partial\Omega}\left(\frac{\partial v}{\partial\nu}w-\frac{\partial w}{\partial\nu}v\right)dS
\ge
\tau^2\int_D\vert v\vert^2dx+O(e^{-2\tau T})+O(\tau e^{-\tau T}e^{-\tau\text{dist}\,(\Omega,B)}).
\end{array}
\tag {2.14}
$$

Here we state a key lemma whose proof is given in Section 4.

\proclaim{\noindent Lemma 2.1.}
It holds that
$$\displaystyle
\liminf_{\tau\longrightarrow\infty}\tau^{6}e^{2\tau\,\text{dist}\,(D,B)}
\int_D\vert v\vert^2dx>0.
\tag {2.15}
$$

\endproclaim

Multiplying the both side of (2.14) by $\tau^{4}e^{2\tau\,\text{dist}\,(D,B)}$, we have
$$\begin{array}{c}
\displaystyle
\tau^{4}e^{2\tau\,\text{dist}\,(D,B)}
\int_{\partial\Omega}\left(\frac{\partial v}{\partial\nu}w-\frac{\partial w}{\partial\nu}v\right)dS
\ge
\tau^{6}e^{2\tau\,\text{dist}\,(D,B)}\int_D\vert v\vert^2dx\\
\\
\displaystyle
+O(\tau^{4}e^{-2\tau(T-\text{dist}\,(D,B))})
+O(\tau^{5}e^{-\tau(T-2\text{dist}\,(D,B)+\text{dist}\,(\Omega,B))})
\end{array}
$$
and thus from (2.15) one gets
$$\displaystyle
\liminf_{\tau\longrightarrow\infty}\tau^{4}e^{2\tau\,\text{dist}\,(D,B)}
\int_{\partial\Omega}\left(\frac{\partial v}{\partial\nu}w-\frac{\partial w}{\partial\nu}v\right)dS>0
\tag {2.16}
$$
provided if $T>2\text{dist}\,(D,B)-\text{dist}\,(\Omega,B)$.

On the other hand, using (2.3), (2.7) and (2.10) one has
$$\begin{array}{c}
\displaystyle
\int_{\partial\Omega}\left(\frac{\partial v}{\partial\nu}w-\frac{\partial w}{\partial\nu}v\right)dS
=\int_D\vert\nabla v\vert^2dx+\tau^2\int_D\vert v\vert^2dx
\\
\\
\displaystyle
+\int_{\partial D}\frac{\partial v}{\partial\nu}(w-v)dS
+O(\tau e^{-\tau T}e^{-\tau\text{dist}\,(\Omega,B)}).
\end{array}
\tag {2.17}
$$

From (2.13) we have, as $\tau\longrightarrow\infty$
$$\displaystyle
\int_D\vert\nabla v\vert^2dx+\tau^2\int_D\vert v\vert^2dx
=O(\tau^2 e^{-2\tau\text{dist}\,(D,B)}).
\tag {2.18}
$$
Concerning with the bound on the third term of the right hand side of (2.17),
we have the following lemma.

\proclaim{\noindent Lemma 2.2.}
It holds that, as $\tau\longrightarrow\infty$
$$\displaystyle
\left\Vert\frac{\partial v}{\partial\nu}\right\Vert_{H^{-1/2}(\partial D)}
=O(\tau^2 e^{-\tau\text{dist}\,(D,B)})
\tag {2.19}
$$
and
$$\displaystyle
\Vert w-v\Vert_{H^1(\Bbb R^3\setminus\overline D)}
=O(\tau e^{-\tau T}+\tau^2 e^{-\tau\text{dist}\,(D,B)}).
\tag {2.20}
$$
\endproclaim
{\it\noindent Proof.}
Set $\epsilon=w-v$.
It follows from (2.9) that
$$\begin{array}{c}
\displaystyle
\int_{\Bbb R^3\setminus\overline D}\vert\nabla\epsilon\vert^2dx+\tau^2\int_{\Bbb R^3\setminus\overline D}\vert\epsilon\vert^2dx
=-e^{-\tau T}\int_{\Bbb R^3\setminus\overline D}(u'(x,T)+\tau u(x,T))\epsilon dx-
\int_{\partial D}\frac{\partial\epsilon}{\partial\nu}\epsilon dS
\end{array}
$$
and thus this yields
$$\begin{array}{c}
\displaystyle
\Vert\nabla\epsilon\Vert_{L^2(\Bbb R^3\setminus\overline D)}^2
+\tau^2\Vert\epsilon\Vert_{L^2(\Bbb R^3\setminus\overline D)}^2
\\
\\
\displaystyle
\le
e^{-\tau T}\Vert u'(\,\cdot\,,T)+\tau u(\,\cdot\,,T)\Vert_{L^2(\Bbb R^3\setminus\overline D)}
\Vert\epsilon\Vert_{L^2(\Bbb R^3\setminus\overline D)}
+\left\vert\int_{\partial D}\frac{\partial v}{\partial\nu}\epsilon dS\right\vert
\end{array}
\tag {2.21}
$$

The boundedness of the trace operator $H^1(\Bbb R^3\setminus\overline D)\longrightarrow H^{1/2}(\partial D)$ yields
$$\begin{array}{c}
\displaystyle
\left\vert\int_{\partial D}\frac{\partial v}{\partial\nu}\epsilon dS\right\vert
\le\left\Vert\frac{\partial v}{\partial\nu}\right\Vert_{H^{-1/2}(\partial D)}
\Vert\epsilon\vert_{\partial D}\Vert_{H^{1/2}(\partial D)}\\
\\
\displaystyle
\le C\left\Vert\frac{\partial v}{\partial\nu}\right\Vert_{H^{-1/2}(\partial D)}
\Vert\epsilon\Vert_{H^{1}(\Bbb R^3\setminus\overline D)},
\end{array}
\tag {2.22}
$$
where $C$ is a positive constant independent of $\tau$.  Now from (2.21) and (2.22) we obtain,
for all $\delta>0$
$$\begin{array}{c}
\displaystyle
(1-C\delta^2)\Vert\nabla\epsilon\Vert_{L^2(\Bbb R^3\setminus\overline D)}^2
+\left(\tau^2-C\delta^2-\frac{\delta^2}{2}\right)\Vert\epsilon\Vert_{L^2(\Bbb R^3\overline D)}^2\\
\\
\displaystyle
\le\frac{\delta^{-2}}{2}e^{-2\tau T}
\Vert u'(\,\cdot\,,T)+\tau u(\,\cdot\,,T)
\Vert_{L^2(\Bbb R^3\setminus\overline D)}^2
+C\delta^{-2}\left\Vert\frac{\partial v}{\partial\nu}\right\Vert_{H^{-1/2}(\partial D)}^2.
\end{array}
\tag {2.23}
$$

Since $v$ satisfies $(\triangle-\tau^2)v=0$ in $D$, for all $\Psi\in H^1(D)$
we have
$$\displaystyle
\int_{\partial D}\frac{\partial v}{\partial\nu}\Psi dS
=\tau^2\int_D v\Psi dx+\int_D\nabla v\cdot\nabla\Psi dx.
$$
The trace operator $H^1(D)\longrightarrow H^{1/2}(\partial D)$ has a bounded right inverse.
This together with the identity above yields
$$\displaystyle
\left\Vert\frac{\partial v}{\partial\nu}\right\Vert_{H^{-1/2}(\partial D)}
\le C(\Vert\nabla v\Vert_{L^2(D)}+\tau^2\Vert v\Vert_{L^2(D)}),
$$
where $C$ is a positive constant independent of $\tau$.  It follows from (2.13)
that this right hand side has the bound $O(\tau^2e^{-\tau \text{dist}(D,B)})$.
Thus (2.19) is valid.   Now (2.20) is a consequence of (2.19) and (2.23).

\noindent
$\Box$

Now we continue the proof of Theorem 1.1.

It follows from (2.19), (2.20) and (2.22) that
$$\displaystyle
\int_{\partial D}\frac{\partial v}{\partial\nu}(w-v)dS
=O(\tau^3 e^{-\tau(T+\text{dist}\,(D,B))}
+\tau^4 e^{-2\tau\text{dist}\,(D,B)}).
$$
From this, (2.17), (2.18) we obtain
$$\begin{array}{c}
\displaystyle
e^{2\tau\text{dist}\,(D,B)}\int_{\partial\Omega}\left(\frac{\partial v}{\partial\nu}w-\frac{\partial w}{\partial\nu}v\right)dS\\
\\
\displaystyle
=
O(\tau^3 e^{-\tau(T-\text{dist}\,(D,B))}
+\tau^4
+\tau e^{-\tau(T-2\text{dist}\,(D,B)+\text{dist}\,(\Omega,B)}).
\end{array}
$$
A combination of this and the estimate $\text{dist}\,(D,B)>\text{dist}\,(\Omega,B)$ gives
$$\displaystyle
\limsup_{\tau\longrightarrow\infty}\tau^{-4}
e^{2\tau\text{dist}\,(D,B)}
\left\vert\int_{\partial\Omega}\left(\frac{\partial v}{\partial\nu}w-\frac{\partial w}{\partial\nu}v\right)
dS\right\vert<\infty
\tag {2.24}
$$
provided $T>2\,\text{dist}\,(D,B)-\text{dist}\,(\Omega,B)$.
Now the formula (1.9) is a direct consequence of (2.16) and (2.24).

\noindent
$\Box$

\section{The enclosure method for penetrable obstacles}

First we specify what we mean by the solution of (1.10).
By Theorem 1 on p.558 in \cite{DL}, given $u^0\in H^1(\Bbb R^3)$
and $u^1\in L^2(\Bbb R^3)$
we know that there exists a unique $u$ satisfying
$$\displaystyle
u\in L^2(0,\,T;H^1(\Bbb R^3)),\,
u'\in L^2(0,\,T;H^1(\Bbb R^3)),\,
u''\in L^2(0,\,T;(H^1(\Bbb R^3))')
$$
such that, for all $\phi\in H^1(\Bbb R^3)$
$$\displaystyle
<u''(t),\phi>
+\int_{\Bbb R^3}\gamma(x)\nabla u(x,t)\cdot\nabla\phi(x)dx=0\,\,\text{a.e.}\,t\in\,]0,\,T[
$$
and $\displaystyle u(x,0)=u^0,\,\,u'(x,0)=u^1$.  In this section we say that this $u$
for $u^0=0$ and $u^1=f$ is the solution of (1.10).

\subsection{A basic identity}

Let $u$ be the solution of (1.10).
Define
$$\displaystyle
w(x;\tau)=\int_0^T e^{-\tau t}u(x,t)dt,\,\,\,x\in\Bbb R^3.
$$
This $w$ belongs to $H^1(\Bbb R^3)$.

From integration by parts (Proposition 2 on p.558 in \cite{DL}) it follows that, for all $\phi\in H^1(\Bbb R^3)$
$$\displaystyle
\int_{\Bbb R^3}\gamma\nabla w\cdot\nabla\phi dx
+\int_{\Bbb R^3}(\tau^2w-f)\phi dx
=-e^{-\tau T}
\int_{\Bbb R^3}(u'(x,T)+\tau u(x,T))\phi dx.
\tag {3.1}
$$

\noindent
This means that in a weak sense $w$ satisfies
$$\displaystyle
(\nabla\cdot\gamma\nabla-\tau^2)w+f(x)=e^{-\tau T}
(u'(x,T)+\tau u(x,T))\,\,\text{in}\,\Bbb R^3.
$$

By a similar reason as the sound-hard obstacle case we know that $w\in H^2_{\text{loc}}(\Bbb R^3\setminus\overline D)$.
Since $\gamma(x)\equiv 1$ in $\Bbb R^3\setminus\overline D$, we define
$\gamma\nabla w\cdot\nu\vert_{\partial\Omega}$ as $\nabla w\vert_{\partial\Omega}\cdot\nu$, where
$\nabla w\vert_{\partial\Omega}$ denotes the trace of $\nabla w$ onto $\partial\Omega$.
Note also that $w$ satisfies $(\triangle-\tau^2)w+f(x)=e^{-\tau T}(u'(x,T)+\tau u(x,T))$ a.e. $x\in\Bbb R^3$.

In this subsection we derive an important identity.

\proclaim{\noindent Proposition 3.1.}
Let $v$ be the weak solution of (1.4).
It holds that
$$\begin{array}{c}
\displaystyle
\int_{\partial\Omega}\{(\nabla v\cdot\nu)w-(\gamma\nabla w\cdot\nu)v\}dS\\
\\
\displaystyle
=-\int_Dh\nabla v\cdot\nabla v dx
+\int_{\Bbb R^3}\gamma\nabla(w-v)\cdot\nabla(w-v)dx
+\tau^2\int_{\Bbb R^3}\vert w-v\vert^2dx-\int_{\Omega}f(w-v)dx\\
\\
\displaystyle
+e^{-\tau T}\int_{\Bbb R^3}(u'(x,T)+\tau u(x,T))(w-v)dx
-e^{-\tau T}\int_{\Omega}(u'(x,T)+\tau u(x,T))vdx.
\end{array}
\tag {3.2}
$$
\endproclaim

{\it\noindent Proof.}
Using a similar argument for the proof of (2.5) and (2.6), we obtain
$$\displaystyle
\int_{\partial\Omega}(\gamma\nabla w\cdot\nu)vdS
=\int_{\Omega}\gamma\nabla w\cdot\nabla vdx+\int_{\Omega}\{\tau^2 w-f+e^{-\tau T}(u'(x,T)+\tau u(x,T)\}vdx
$$
and
$$\displaystyle
\int_{\partial\Omega}(\nabla v\cdot\nu)wdS
=\int_{\Omega}\nabla v\cdot\nabla w dx
+\int_{\Omega}(\tau^2 v-f)wdx.
$$
From these we obtain
$$\begin{array}{c}
\displaystyle
\int_{\partial\Omega}\{(\nabla v\cdot\nu)w-(\gamma\nabla w\cdot\nu)v\}dS\\
\\
\displaystyle
=-\int_{D}h\nabla w\cdot\nabla vdx
-\int_{\Omega}f(w-v)dx
-e^{-\tau T}\int_{\Omega}(u'(x,T)+\tau u(x,T))vdx.
\end{array}
\tag {3.3}
$$

Write
$$\displaystyle
-\int_D h\nabla w\cdot\nabla vdx
=-\int_D h\nabla v\cdot\nabla vdx
-\int_Dh\nabla(w-v)\cdot\nabla vdx.
\tag {3.4}
$$
Since $v$ satisfies (1.5),
from (3.1) we have, for all $\phi\in H^1(\Bbb R^3)$
$$\begin{array}{c}
\displaystyle
-\int_{\Bbb R^3}\gamma\nabla(w-v)\cdot\nabla\phi dx
-\tau^2\int_{\Bbb R^3}(w-v)\phi dx\\
\\
\displaystyle
=-\int_{\Bbb R^3}(I_3-\gamma)\nabla v\cdot\nabla\phi dx
+e^{-\tau T}\int_{\Bbb R^3}(u'(x,T)+\tau u(x,T))\phi dx.
\end{array}
\tag {3.5}
$$
This means that the $w-v$ satisfies, in a weak sense
$$
\displaystyle
(\nabla\cdot\gamma\nabla-\tau^2)(w-v)
=\nabla\cdot(I_3-\gamma)\nabla v+e^{-\tau T}(u'(x,T)+\tau u(x,T))
\,\,\text{in}\,\Bbb R^3.
$$
Substituting $w-v$ for $\phi$ in (3.5), we obtain
$$\begin{array}{l}
\displaystyle
\int_{\Bbb R^3}(I_3-\gamma)\nabla v\cdot\nabla(w-v)dx
=\int_{\Bbb R^3}\gamma\nabla(w-v)\cdot\nabla(w-v)dx\\
\\
\displaystyle
+\tau^2\int_{\Bbb R^3}(w-v)(w-v)dx
+e^{-\tau T}\int_{\Bbb R^3}(u'(x,T)+\tau u(x,T))(w-v)dx.
\end{array}
\tag {3.6}
$$
A combination of (3.4) and (3.6) gives
$$\begin{array}{l}
\displaystyle
-\int_D h\nabla w\cdot\nabla vdx
=-\int_D h\nabla v\cdot\nabla vdx\\
\\
\displaystyle
+\int_{\Bbb R^3}\gamma\nabla(w-v)\cdot\nabla(w-v)dx
+\tau^2\int_{\Bbb R^3}(w-v)(w-v)dx\\
\\
\displaystyle
+e^{-\tau T}\int_{\Bbb R^3}(u'(x,T)+\tau u(x,T))(w-v)dx.
\end{array}
$$
Now from this and (3.3) we obtain (3.2).

\noindent
$\Box$

In particular, choose $f$ in such a way that $\text{supp}\,f\cap\overline\Omega=\emptyset$.
Then (3.2) become
$$\begin{array}{c}
\displaystyle
\int_{\partial\Omega}
\{(\nabla v\cdot\nu)w-(\gamma\nabla w\cdot\nu)v\}dS\\
\\
\displaystyle
=-\int_D h\nabla v\cdot\nabla v dx+\int_{\Bbb R^3}\gamma\nabla(w-v)\cdot\nabla(w-v)dx
+\tau^2\int_{\Bbb R^3}\vert w-v\vert^2dx\\
\\
\displaystyle
+e^{-\tau T}\int_{\Bbb R^3}(u'(x,T)+\tau u(x,T))(w-v)dx
-e^{-\tau T}
\int_{\Omega}(u'(x,T)+\tau u(x,T))vdx.
\end{array}
\tag {3.7}
$$
This is our first basic identity which is useful in the proof of Theorem 1.2 under the assumption (A.1).

Unfortunately, for (A2) this identity does not work.  However, one can rewrite this by replacing the role of
$v$ and $w$ in the proof of Proposition 3.1.  
More precisely, set $\tilde{f}=f-e^{-\tau T}(u'(x,T)+\tau u(x,T))$.
The points are: $w$ satisfies 
$\displaystyle\nabla\cdot\gamma\nabla w-\tau^2 w+\tilde{f}=0\,\,\text{in}\,\Bbb R^3$
and $v$ satisfies $\displaystyle
\nabla\cdot I_3\nabla v-\tau^2 v+\tilde{f}
=-e^{-\tau T}(u'(x,T)+\tau u(x,T))\,\,\text{in}\,\Bbb R^3$.
Thus changing the role of $v$ and $w$ in the proof of Proposition 3.1, we can easily obtain another expression of (3.2).
\proclaim{\noindent Proposition 3.2.}
Let $v$ be the weak solution of (1.4).
It holds that
$$\begin{array}{c}
\displaystyle
\int_{\partial\Omega}\{(\gamma\nabla w\cdot\nu)v-(\nabla v\cdot\nu)w\}dS\\
\\
\displaystyle
=\int_Dh\nabla w\cdot\nabla wdx
+\int_{\Bbb R^3}\nabla(v-w)\cdot\nabla(v-w)dx
+\tau^2\int_{\Bbb R^3}\vert v-w\vert^2dx
-\int_{\Omega}f(v-w)dx\\
\\
\displaystyle
-e^{-\tau T}\int_{\Bbb R^3}(u'(x,T)+\tau u(x,T))(v-w)dx
+e^{-\tau T}\int_{\Omega}(u'(x,T)+\tau u(x,T))vdx.
\end{array}
\tag {3.8}
$$
\endproclaim

In particular, if $\text{supp}\,f\cap\overline\Omega=\emptyset$, then (3.8) gives
$$\begin{array}{c}
\displaystyle
\int_{\partial\Omega}\{(\gamma\nabla w\cdot\nu)v-(\nabla v\cdot\nu)w\}dS\\
\\
\displaystyle
=\int_Dh\nabla w\cdot\nabla wdx
+\int_{\Bbb R^3}\nabla(v-w)\cdot\nabla(v-w)dx
+\tau^2\int_{\Bbb R^3}\vert v-w\vert^2dx\\
\\
\displaystyle
-e^{-\tau T}\int_{\Bbb R^3}(u'(x,T)+\tau u(x,T))(v-w)dx
+e^{-\tau T}\int_{\Omega}(u'(x,T)+\tau u(x,T))vdx.
\end{array}
\tag {3.9}
$$

\subsection{Proof of Theorem 1.2.}
First we consider the case when (A1) is satisfied.  Using (A1) and the identity
$$\begin{array}{c}
\displaystyle
\tau^2\vert w-v\vert^2+e^{-\tau t}(w-v)(u'(x,T)+\tau u(x,T))\\
\\
\displaystyle
=\left\vert\tau(w-v)+\frac{e^{-\tau T}}{2\tau}(u'(x,T)+\tau u(x,T))\right\vert^2
-\frac{e^{-2\tau T}}{4\tau^2}\vert u'(x,T)+\tau u(x,T)\vert^2,
\end{array}
$$
we have from (3.7)
$$\begin{array}{c}
\displaystyle
\int_{\partial\Omega}
\{(\nabla v\cdot\nu)w-(\gamma\nabla w\cdot\nu)v\}dS
\ge C\int_D\vert\nabla v\vert^2dx
\\
\\
\displaystyle
-\frac{e^{-2\tau T}}{4\tau^2}\int_{\Bbb R^3}\vert u'(x,T)+\tau u(x,T)\vert^2dx
-e^{-\tau T}\int_{\Omega}
(u'(x,T)+\tau u(x,T))vdx.
\end{array}
\tag {3.10}
$$
Since we have
$$\displaystyle
\int_{\Bbb R^3}\vert u'(x,T)+\tau u(x,T)\vert^2dx=O(\tau^2)
$$
and the estimate $\Vert v\Vert_{L^2(\Omega)}
=O(e^{-\tau\text{dist}\,(\Omega,B)})$, it follows from (3.10) that
$$\begin{array}{c}
\displaystyle
\int_{\partial\Omega}
\{(\nabla v\cdot\nu)w-(\gamma\nabla w\cdot\nu)v\}dS
\ge
C\int_D\vert\nabla v\vert^2dx+O(e^{-2\tau T})+O(\tau e^{-\tau T}e^{-\tau\text{dist}\,(\Omega,B)}).
\end{array}
\tag {3.11}
$$

Here we state a key lemma whose proof is given in the next section.

\proclaim{\noindent Lemma 3.1.}
It holds that
$$\displaystyle
\liminf_{\tau\longrightarrow\infty}\tau^{4}e^{2\tau\,\text{dist}\,(D,B)}
\int_D\vert\nabla v\vert^2dx>0.
\tag {3.12}
$$

\endproclaim

Multiplying the both side of (3.11) by $\tau^{\mu}e^{2\tau\,\text{dist}\,(D,B)}$, we have
$$\begin{array}{c}
\displaystyle
\tau^{4}e^{2\tau\,\text{dist}\,(D,B)}
\int_{\partial\Omega}
\{(\nabla v\cdot\nu)w-(\gamma\nabla w\cdot\nu)v\}dS
\ge
\tau^{4}e^{2\tau\,\text{dist}\,(D,B)}\int_D\vert\nabla v\vert^2dx\\
\\
\displaystyle
+O(\tau^{4}e^{-2\tau(T-\text{dist}\,(D,B))})
+O(\tau^{5}e^{-\tau(T-2\text{dist}\,(D,B)+\text{dist}\,(\Omega,B))})
\end{array}
$$
and thus from (3.12) one gets
$$\displaystyle
\liminf_{\tau\longrightarrow\infty}\tau^{4}e^{2\tau\,\text{dist}\,(D,B)}
\int_{\partial\Omega}
\{(\nabla v\cdot\nu)w-(\gamma\nabla w\cdot\nu)v\}dS
>0
\tag {3.13}
$$
provided if $T>2\text{dist}\,(D,B)-\text{dist}\,(\Omega,B)$.

On the other hand, using (3.3) and (3.4), one gets
$$\begin{array}{c}
\displaystyle
\int_{\partial\Omega}
\{(\nabla v\cdot\nu)w-(\gamma\nabla w\cdot\nu)v\}dS\\
\\
\displaystyle
=-\int_D h\nabla v\cdot\nabla vdx-\int_Dh\nabla(w-v)\cdot\nabla vdx
+O(\tau e^{-\tau T}e^{-\tau\text{dist}\,(\Omega,B)}).
\end{array}
\tag {3.14}
$$

From (2.13) we have, as $\tau\longrightarrow\infty$
$$\displaystyle
\int_D\vert\nabla v\vert^2dx
=O(\tau^2 e^{-2\tau\text{dist}\,(D,B)}).
\tag {3.15}
$$
Concerning with the bound on the second term of the right hand side of (3.12),
we have the following lemma.

\proclaim{\noindent Lemma 3.2.}
It holds that, as $\tau\longrightarrow\infty$
$$\displaystyle
\Vert w-v\Vert_{H^1(\Bbb R^3)}
=O(\tau e^{-\tau T}+\tau e^{-\tau\text{dist}\,(D,B)})
\tag {3.16}
$$
\endproclaim
{\it\noindent Proof.}
Set $\epsilon=w-v$.  From (3.6) we have
$$\displaystyle
-\int_{D}h\nabla v\cdot\nabla\epsilon dx
-\int_{\Bbb R^3}\gamma\nabla\epsilon\cdot\nabla\epsilon dx
-\tau^2\int_{\Bbb R^3}\vert\epsilon\vert^2dx
=e^{-\tau T}\int_{\Bbb R^3}(u'(x,T)+\tau u(x,T))\epsilon dx.
$$
This yields
$$\begin{array}{c}
\displaystyle
C\Vert\nabla\epsilon\Vert_{L^2(\Bbb R^3)}^2
+\tau^2\Vert\epsilon\Vert_{L^2(\Bbb R^3)}^2
\\
\\
\displaystyle
\le
e^{-\tau T}\Vert u'(\,\cdot\,,T)+\tau u(\,\cdot\,,T)\Vert_{L^2(\Bbb R^3)}
\Vert\epsilon\Vert_{L^2(\Bbb R^3)}
+\left\vert\int_{D}h\nabla v\cdot\nabla\epsilon dx\right\vert
\end{array}
$$
and thus one gets,
for all $\delta>0$
$$\begin{array}{c}
\displaystyle
\left(C-\frac{\delta^2}{2}\right)\Vert\nabla\epsilon\Vert_{L^2(\Bbb R^3)}^2
+\left(\tau^2-\frac{\delta^2}{2}\right)\Vert\epsilon\Vert_{L^2(\Bbb R^3)}^2\\
\\
\displaystyle
\le\frac{\delta^{-2}}{2}e^{-2\tau T}
\Vert u'(\,\cdot\,,T)+\tau u(\,\cdot\,,T)
\Vert_{L^2(\Bbb R^3)}^2
+\frac{\delta^{-2}}{2}\Vert h\Vert_{L^{\infty}(D)}^2\Vert\nabla v\Vert_{L^2(D)}^2.
\end{array}
$$
Now a combination of this and (3.15) yields (3.16).

\noindent
$\Box$

A combination of (3.15) and (3.16) gives
$$\displaystyle
\int_Dh\nabla(w-v)\cdot\nabla vdx
=O(\tau^2 e^{-\tau(T+\text{dist}\,(D,B))}
+\tau^2 e^{-2\tau\text{dist}\,(D,B)}).
$$
From this, (3.14), (3.15) we obtain
$$\begin{array}{c}
\displaystyle
e^{2\tau\text{dist}\,(D,B)}\int_{\partial\Omega}\left((\nabla v\cdot\nu)w-
(\gamma\nabla w\cdot\nu)v\right)dS\\
\\
\displaystyle
=
O(\tau^2 e^{-\tau(T-\text{dist}\,(D,B))}
+\tau^2
+\tau e^{-\tau(T-2\text{dist}\,(D,B)+\text{dist}\,(\Omega,B)})).
\end{array}
$$
This together with the estimate $\text{dist}\,(D,B)>\text{dist}\,(\Omega,B)$
yields
$$\displaystyle
\limsup_{\tau\longrightarrow\infty}\tau^{-2}
e^{2\tau\text{dist}\,(D,B)}
\left\vert\int_{\partial\Omega}\left((\nabla v\cdot\nu)w-(\gamma\nabla w\cdot\nu)v\right)
dS\right\vert<\infty
\tag {3.17}
$$
provided $T>2\,\text{dist}\,(D,B)-\text{dist}\,(\Omega,B)$.
Now the conclusion of Theorem 1.2 is a direct consequence of (3.13)
and (3.17).

\noindent
$\Box$

Finally we give a comment on the case when (A2) is satisfied.
In this case we make use of (3.9) instead of (3.7).
A combination of the well known inequality (see \cite{IS})
$$\displaystyle
(\gamma(x)-I_3)\nabla w\cdot\nabla w+\nabla(w-v)\cdot\nabla(w-v)
\ge(\gamma(x)-I_3)\gamma(x)^{-1/2}\nabla v\cdot\gamma(x)^{-1/2}\nabla v
$$
and (A2) yields that there exists a positive constant $C$ such that 
$$\displaystyle
h(x)\nabla w\cdot\nabla w+\nabla(w-v)\cdot\nabla(w-v)\ge C\vert\nabla v\vert^2
$$
for a.e. $x\in D$.  This together with (3.9) gives the lower estimate
$$\begin{array}{c}
\displaystyle
\int_{\partial\Omega}
\{(\gamma\nabla w\cdot\nu)v-(\nabla v\cdot\nu)w\}dS
\ge C\int_D\vert\nabla v\vert^2dx\\
\\
\displaystyle
-\frac{e^{-2\tau T}}{4\tau^2}\int_{\Bbb R^3}\vert u'(x,T)+\tau u(x,T)\vert^2dx
+e^{-\tau T}\int_{\Omega}
(u'(x,T)+\tau u(x,T))vdx.
\end{array}
$$
which corresponds to (3.10).  Applying the arugument for the proof of (3.13) to this right hand side
we obain
$$\displaystyle
\liminf_{\tau\longrightarrow\infty}\tau^{4}e^{2\tau\,\text{dist}\,(D,B)}
\int_{\partial\Omega}
\{(\gamma\nabla w\cdot\nu)v-(\nabla v\cdot\nu)w\}dS
>0.
$$
Since (3.17) is valid also for case (A2) we obtain the desired conclusion.

\noindent
$\Box$

\section{Proof of Lemmas 2.1 and 3.1}

\subsection{Proof of Lemma 2.1.}
Choose points $x_0\in\partial D$ and $y_0\in\partial B$ such that
$\text{dist}\,(D,B)=\vert x_0-y_0\vert$.  Since we have assumed that $\partial D$ is smooth,
one can find an open ball $B'$ such that $B'\subset D$ and $x_0\in\partial B'\cap\partial D$.
Since $\text{dist}\,(B',B)=\vert x_0-y_0\vert$, it suffices to prove (2.15) in the case when $D=B'$.

Write
$$\displaystyle
v(x)^2
=\left(\frac{1}{4\pi}\right)^2
\int_{B\times B}
\frac{e^{-\tau(\vert y_1-x\vert+\vert y_2-x\vert)}}{\vert x-y_1\vert \vert x-y_2\vert}f(y_1)f(y_2)dy_1dy_2.
$$
It follows from the assumption on $f$ that
$$\displaystyle
v(x)^2\ge C^2I(x,\tau)^2
\tag {4.1}
$$
where
$$\displaystyle
I(x,\tau)
=\frac{1}{4\pi}\int_B \frac{e^{-\tau\vert y-x\vert}}{\vert y-x\vert}dy,\,x\in B'.
$$

We denote by $p$ and $\eta$ the center and radius of $B$, respectively.
Using the polar coordinates centered at $x$, one can write
$$\displaystyle
B=\{y=x+r\omega\,\vert\,\omega\in S(x,B),\,r^{+}(\omega)<r<r^{-}(\omega)\}
$$
where
$$\displaystyle
S(x,B)=\{\omega\in S^2\,\vert\,\omega\cdot(p-x)>\sqrt{\vert x-p\vert^2-\eta^2}\}
$$
and
$$\displaystyle
r^{\pm}(\omega)=\omega\cdot(p-x)\mp\sqrt{(\omega\cdot(p-x))^2-\vert x-p\vert^2+\eta^2}.
$$
This together with $r^{+}(\omega)\ge\text{dist}\,(B',B)$ yields
$$\displaystyle
I(x,\tau)
\ge \text{dist}\,(B',B)\int_{S(x,B)}d\omega
\int_{r^{+}(\omega)}^{r^{-}(\omega)}e^{-\tau r}dr
$$
and thus we have
$$\displaystyle
\tau I(x,\tau)
\ge\frac{\text{dist}\,(B',B)}{4\pi}
\left(\int_{S(x,B)}e^{-\tau r^+(\omega)}d\omega-\int_{S(x,B)}e^{-\tau r^-(\omega)}d\omega\right).
\tag {4.2}
$$

Define $d_B(x)=\inf\{\vert x-y\vert\,\vert\,y\in\,B\}$.
First we give an estimate for the second integral in the right hand side of (4.2).
Since $d_B(x)=\vert x-p\vert-\eta$ and $r^{-}(\omega)>\sqrt{\vert x-p\vert^2-\eta^2}$, we have
$$\begin{array}{c}
\displaystyle
\int_{S(x,B)}e^{-\tau r^{-}(\omega)}d\omega
\le\int_{S(x,B)}e^{-\tau\sqrt{d_B(x)}\sqrt{\vert x-p\vert+\eta}}d\omega\\
\\
\displaystyle
\le 4\pi e^{-\tau d_B(x)}
e^{-\tau(\sqrt{d_B(x)}\sqrt{\vert x-p\vert+\eta}-d_B(x))}.
\end{array}
$$
Here note that
$$\displaystyle
\sqrt{d_B(x)}\sqrt{\vert x-p\vert+\eta}-d_B(x)
=\frac{2\eta\sqrt{d_B(x)}}
{\sqrt{\vert x-p\vert+\eta}+\sqrt{d_B(x)}}>0.
$$
Thus we obtain
$$\displaystyle
e^{\tau d_B(x)}\int_{S(x,B)}e^{-\tau r^{-}(\omega)}d\omega
\le 4\pi e^{-A\tau}
\tag {4.3}
$$
where
$$\displaystyle
A=\inf_{x\in B'}\frac{2\eta\sqrt{d_B(x)}}
{\sqrt{\vert x-p\vert+\eta}+\sqrt{d_B(x)}}.
$$
The points are: $A$ is positive and independent of $x\in B'$.

Next we consider the first integral of the right hand side of (4.2).
The surface $S(x,B)$ has the parameterization:
$$\displaystyle
\omega(r,\theta)
=r\cos\,\theta\mbox{\boldmath $a$}(x)+r\sin\,\theta\,\mbox{\boldmath $b$}(x)
+\frac{p-x}{\vert p-x\vert}\sqrt{1-r^2},
\tag {4.4}
$$
where $0<r<\eta/\vert x-p\vert$ and $0\le \theta<2\pi$; $\mbox{\boldmath $a$}(x)$ and $\mbox{\boldmath $b$}(x)$
are unit vectors perpendicular each other and satisfy $\mbox{\boldmath $a$}(x)\times\mbox{\boldmath $b$}(x)
=(p-x)/\vert p-x\vert$.
Since $\displaystyle d\omega=(r/\sqrt{1-r^2})drd\theta$, one can write
$$\displaystyle
\int_{S(x,B)}e^{-\tau r^+(\omega)}d\omega
=\int_0^{2\pi}d\theta\int_0^{\eta/\vert x-p\vert}
e^{-\tau r^+(\omega(r,\theta))}\frac{rdr}{\sqrt{1-r^2}}.
$$
Here we note that $r^{+}(\omega(r,\theta))$
is independent of $\theta$.  In fact we have
$$\displaystyle
r^+(\omega(r,\theta))=\vert p-x\vert\sqrt{1-r^2}-\sqrt{\eta^2-r^2\vert p-x\vert^2}.
$$
Thus this yields
$$\displaystyle
\int_{S(x,B)}e^{-\tau r^+(\omega)}d\omega
\ge 2\pi\int_0^{\eta/\vert p-x\vert}
e^{-\tau r^{+}(\omega(r,0))}rdr.
\tag {4.5}
$$

Since $d_B(x)=\vert p-x\vert-\eta$, we obtain
$$\begin{array}{c}
\displaystyle
r^+(\omega(r,0))-d_B(x)
=\vert p-x\vert\sqrt{1-r^2}-\sqrt{\eta^2-r^2\vert p-x\vert^2}
-\vert p-x\vert+\eta\\
\\
\displaystyle
=\vert p-x\vert(\sqrt{1-r^2}-1)
+(\eta-\sqrt{\eta^2-r^2\vert p-x\vert^2})\\
\\
\displaystyle
=r^2\vert p-x\vert\left(\frac{1}{\eta+\sqrt{\eta^2-r^2\vert p-x\vert^2}}-\frac{1}{1+\sqrt{1-r^2}}\right)\\
\\
\displaystyle
\le\frac{r^2\vert p-x\vert}{\eta+\sqrt{\eta^2-r^2\vert p-x\vert^2}}
<\frac{r^2\vert p-x\vert}{\eta}.
\end{array}
\tag {4.6}
$$
Now set $\displaystyle M=\sup_{x\in B'}(\vert p-x\vert/\eta)$.
Since $r^+(\omega(r,0))-d_B(x)\ge 0$, from (4.6) we obtain
$$\begin{array}{c}
\displaystyle
\int_0^{\eta/\vert x-p\vert}
e^{-\tau r^+(\omega(r,0))}rdr
\ge e^{-\tau d_B(x)}
\int_0^{\eta/\vert p-x\vert}
e^{-\tau Br^2}rdr\\
\\
\displaystyle
=\frac{e^{-\tau d_B(x)}}{2\tau M}(1-e^{-\tau B(\eta/\vert p-x\vert)^2})\\
\\
\displaystyle
\ge
\frac{e^{-\tau d_B(x)}}{2\tau M}(1-e^{-\tau/B}).
\end{array}
\tag {4.7}
$$
From this, (4.2), (4.3) and (4.5) we can conclude that: there exist $\tau_0>0$ and $C'>0$
independent of $x\in B'$ such that, for all $\tau\ge\tau_0$
$$\displaystyle
\tau^2 e^{\tau d_B(x)}I(x,\tau)\ge C'.
$$
This together with (4.1) yields
$$\displaystyle
\tau^4\int_{B'}v(x)^2dx
\ge (CC')^{2}\int_{B'}e^{-2\tau d_B(x)}dx
=(CC')^{2}e^{2\tau\eta}\int_{B'}e^{-2\tau\vert x-p\vert}dx.
$$
We have already known (Proposition 3.2 in \cite{IK2}) that
$$\displaystyle
\liminf_{\tau\longrightarrow\infty}\tau^2 e^{2\tau d_{B'}(p)}\int_{B'}e^{-2\tau\vert x-p\vert}dx>0.
\tag {4.8}
$$
This yields
$$\displaystyle
\liminf_{\tau\longrightarrow\infty}\tau^6e^{2\tau(d_{B'}(p)-\eta)}\int_{B'}v(x)^2dx>0.
$$
Since $d_{B'}(p)-\eta=\text{dist}\,(B',B)$ we conclude that
(2.15) is valid.

\noindent
$\Box$

\subsection{Proof of Lemma 3.1.}

We employ the same notation used in the proof of Lemma 2.1.
It suffices to prove (3.12) in the case when $D=B'$.
Since
$$\displaystyle
\nabla v(x)
=-\frac{1}{4\pi}
\int_{B}
\frac{1+\tau\vert x-y\vert^2}{\vert x-y\vert^3}
e^{-\tau\vert x-y\vert}f(y)(x-y)dy,
$$
we have
$$\begin{array}{c}
\displaystyle
(4\pi)^2\vert\nabla v(x)\vert^2\\
\\
\displaystyle
=
\int_{B\times B}
K(x,y_1)K(x,y_2)e^{-\tau(\vert x-y_1\vert+\vert x-y_2\vert)}
f(y_1)f(y_2)(x-y_1)\cdot(x-y_2)dy_1dy_2,
\end{array}
\tag {4.9}
$$
where
$$\displaystyle
K(x,y)
=\frac{1+\tau\vert x-y\vert^2}
{\vert x-y\vert^3}.
$$
Define
$$\displaystyle
J(x,\tau)
=\int_{B}
\frac{(x-y)}{\vert x-y\vert^3}
e^{-\tau\vert x-y\vert}dy.
$$
Since it holds that
$(x-y_1)\cdot(x-y_2)>0$ for all $y_1,y_2\in B$ and $x\in\Bbb R^3\setminus\overline B$
and $f(y_1)f(y_2)\ge C^2$ a.e. $y_1, y_2\in B$,
from (4.8) we obtain
$$\begin{array}{c}
\displaystyle
(4\pi)^2\vert\nabla v(x)\vert^2\\
\\
\displaystyle
\ge
C^2
\int_{B\times B}
K(x,y_1)K(x,y_2)e^{-\tau(\vert x-y_1\vert+\vert x-y_2\vert)}(x-y_1)\cdot(x-y_2)dy_1dy_2\\
\\
\displaystyle
\ge
C^2(1+\tau\text{dist}\,(D,B)^2)^2
\int_{B\times B}
\frac{(x-y_1)\cdot(x-y_2)}{\vert x-y_1\vert^3\vert x-y_2\vert^3}
e^{-\tau(\vert x-y_1\vert+\vert x-y_2\vert)}dy_1dy_2\\
\\
\displaystyle
=C^2(1+\tau\text{dist}\,(D,B)^2)^2\vert J(x,\tau)\vert^2.
\end{array}
\tag {4.10}
$$
One can write
$$\begin{array}{c}
\displaystyle
J(x,\tau)
=\int_{S(x,B)}\omega d\omega\int_{r^+(\omega)}^{r^{-}(\omega)}e^{-\tau r}dr\\
\\
\displaystyle
=\frac{1}{\tau}
\left(
\int_{S(x,B)}e^{-\tau r^{+}(\omega)}\omega d\omega
-\int_{S(x,B)}e^{-\tau r^{-}(\omega)}\omega d\omega\right).
\end{array}
\tag {4.11}
$$
Since $r^{+}(\omega)\ge d_B(x)$, we have
$$\displaystyle
\left\vert\int_{S(x,B)}e^{-\tau r^{+}(\omega)}\omega d\omega\right\vert
\le 4\pi e^{-\tau d_B(x)}.
$$
The second integral in the right hand side of (4.11) has the bound
$O(e^{-\tau d_B(x)}e^{-A\tau})$.
These give
$$\displaystyle
\tau^2\vert J(x,\tau)\vert^2
=\left\vert\int_{S(x,B)}e^{-\tau r^{+}(\omega)}\omega d\omega\right\vert^2
+O(e^{-A\tau}e^{-2\tau d_B(x)}).
\tag {4.12}
$$
Using parameterization (4.4) of $S(x,B)$, one has
$$\displaystyle
\int_{S(x,B)}e^{-\tau r^+(\omega)}\omega d\omega
=\int_0^{2\pi}d\theta\int_0^{\eta/\vert x-p\vert}
e^{-\tau r^+(\omega(r,\theta))}\omega(r,\theta)\frac{rdr}{\sqrt{1-r^2}}.
$$
Since $r^{+}(\omega(r,\theta))$ is independent of $\theta$ and
$$\displaystyle
\int_0^{2\pi}\omega(r,\theta)d\theta=\frac{p-x}{\vert p-x\vert}\sqrt{1-r^2}2\pi,
$$
we obtain
$$\displaystyle
\int_{S(x,B)}e^{-\tau r^+(\omega)}\omega d\omega
=2\pi\frac{p-x}{\vert p-x\vert}\int_0^{\eta/\vert p-x\vert}
e^{-\tau r^{+}(\omega(r,0))}rdr.
$$
From this together with (4.7) and (4.12)
one can conclude that: there exist $\tau_0>0$ and $C'>0$
independent of $x\in B'$ such that, for all $\tau\ge\tau_0$
$$\displaystyle
\tau^4 e^{\tau 2d_B(x)}\vert J(x,\tau)\vert^2\ge C'.
$$
Thus from (4.10) we obtain
$$\displaystyle
\tau^2\int_{B'}\vert\nabla v(x)\vert^2dx
\ge C''\int_{B'}e^{-2\tau d_B(x)}dx
=C''e^{2\tau\eta}\int_{B'}e^{-2\tau\vert x-p\vert}dx,
$$
where $C''$ is a positive constant.
Hereafter using (4.8), we obtain (3.12).

\noindent
$\Box$

\section{Conclusion and further problems}

In this paper we introduce a simple method for some class of inverse obstacle scattering problems
that employs the values of the wave field over a {\it finite} time interval on a known surface surrounding unknown obstacles
as the observation data.  The wave field is generated by an initial data localized outside the surface
and its form is not specified except for the condition on the support.  The method yields information
about the location and shape of the obstacles more than the convex hull.

\noindent $\bullet$  It would be interesting to apply the method
presented in this paper to other time dependent problems in
electromagnetism, linear elasticity, classical fluids etc.. Those
applications belong to our future plan.

$$\quad$$

\centerline{{\bf Acknowledgements}}

This research was partially supported by Grant-in-Aid for
Scientific Research (C)(No. 21540162) of Japan  Society for the
Promotion of Science.

\end{document}